# Robust Control of General Linear Delay Systems under Dissipativity: Part I – A KSD-based Framework

Qian Feng[1], Wei Xing Zheng[2], Xiaoyu Wang[3] and Feng Xiao[1]

***Abstract*— This paper introduces an effective framework for designing memoryless dissipative full-state feedback for general linear delay systems via the Krasovskiĭ functional (KF) approach, where an arbitrary finite number of pointwise and general distributed delays (DDs) exists in the state, input and output. To handle the infinite dimensionality of DDs, we employ the Kronecker-Seuret Decomposition (KSD) which we recently proposed for analyzing matrix-valued functions in the context of delay systems. The KSD enables factorization or least-squares approximation of any number of $\mathscr{L}^2$ DD kernels from any number of DDs without introducing conservatism. This also facilitates the construction of a complete-type KF with flexible integral kernels by means of a novel integral inequality derived from the least-squares principle. Our solution includes two theorems and an iterative algorithm to compute controller gains without relying on nonlinear solvers. A numerical example is tested to show the effectiveness of the proposed approach.**

## I. Introduction

Pointwise and distributed delays (DDs) are frequently employed to model transport, propagation and aftereffects in a dynamical system. The nature of a pointwise delay is elucidated in [1] as a transport equation coupled with boundary conditions. Meanwhile, delays can also arise from transporting media with more complex structures [2]. A DD is denoted by an integral $\int_{-r}^{0} F(\tau)\boldsymbol{x}(t+\tau)\mathrm{d}\tau$ over a delay interval $[-r, 0]$ with a matrix-valued function $F(\cdot)$, which takes into account a segment of the past dynamics. Systems with both pointwise and DDs have diverse applications such as synchronization of complex networks [3], neural networks [4], and modeling event-triggered mechanisms [5].

Most methods for linear time-delay systems (LTDSs) are carried out in the time or frequency domain, using real or complex analysis. For time-domain methods [6], [7], the Krasovskiĭ functional (KF) approach has been shown to be effective for the stability analysis and stabilization of LTDSs [8]–[12]. This approach seeks to convert the original problems into solving convex semidefinite programming (SDP) problems, which can be computed by efficient numerical algorithms [13]. For a comprehensive collection of the existing literature on this subject, the reader is referred to [6], [14]. In contrast to the Lyapunov approach for an LTI system, the KF approach could only establish sufficient conditions, where the induced conservatism mainly depends on the conservatism of the predetermined form of KFs [7] and the integral inequalities [9] utilized to construct these sufficient conditions. As more general KFs [8], [10] are increasingly adopted to reduce conservatism, congruent transformations may not be applicable in formulating convex controller/observer synthesis conditions from the original stability analysis condition. Finally, a notable method combining both time and frequency domain approaches has been proposed in [15] for the stabilization of LTDSs with DDs, based on an application of the concept of smoothed spectral abscissa [16] and delay Lyapunov matrix [17].

Nevertheless, it is fair to say that there are no effective solutions in the literature for the control of LTDSs with an arbitrary finite number of pointwise and general DDs. Even when considering stability analysis alone, most existing KF approaches impose conservative constraints on the mathematical structure of the state space matrices [18], [19] or limit the number of DD kernels [9] or delays [8], [10]. The method in [15] requires the computation of delay Lyapunov matrix and its derivatives, but the authors did not elaborate on how this computation can be carried out for an LTDS with general DDs or non-commensurate delays. Finally, the solution to the infinite-horizon linear-quadratic control [20], [21] of LTDSs can be obtained by solving operator Riccati equations using the $C_0$ semigroup theory. However, solutions to these equations cannot be explicitly computed and require using sophisticated finite-dimensional approximations[a] [23] to obtain approximate results by solving finite-dimensional algebraic Riccati equations [24], [25]. In fact, special conditions [26, (A1)–(A3)] on the abstract operators of the semigroup representation must be enforced to ensure that these approximation schemes have strong convergence properties, so that the solution to the finite-dimensional algebraic Riccati equations approaches the solution to the operator Riccati equations in norm [27, Section 4.2].

In this paper, we introduce a framework inspired by the KSD concept we recently proposed in [28] in solving general estimator design problem, to compute dissipative controllers for LTDSs with general delay-structures by constructing a complete type KF. The generality of the system model is attained by incorporating an arbitrary finite number of pointwise and general DDs at the states, inputs, and

This work was partially supported by the National Natural Science Foundation of China under Grant Nos. 62303180, 62273145 and 62373150, and Fundamental Research Funds for Central Universities under Grant 2025JC009, China.

1. School of Control and Computer Engineering, North China Electric Power University, Beijing, China. Email: qianfeng@ncepu.edu.cn, qfen204@aucklanduni.ac.nz, fengxiao@ncepu.edu.cn

2. School of Computer, Data and Mathematical Sciences, Western Sydney University, Sydney, Australia. W.Zheng@westernsydney.edu.au

3. Institute of Systems Science, Beijing Wuzi University, Beijing, China. wangxiaoyu@bwu.edu.cn

[a] Galerkin-type approximations (averaging, splines, orthogonal functions, eigenfunctions)

outputs, where the DDs can contain any number of $\mathscr{L}^2$ functions over bounded intervals. Furthermore, we employ the Carathéodory interpretation [29, section 2.6] for our systems differential equations with respect to (w.r.t) the Lebesgue measure, which is better suited for modeling the dynamics of engineering systems often subject to noise and glitches. The use of KSD allows for the decomposition of integral kernel matrices of any DD containing $\mathscr{L}^2$ functions as products of constant matrices and a list of basis functions with specific properties. Benefiting from the structures of KSD, we can construct KFs with general integral kernels in $\mathscr{W}^{1,2}$ as long as they are linearly independent. Once our KF is successfully constructed using integral inequalities derived from the least-squares principle [30, page 182], a controller synthesis condition is formulated for the dissipative control problem via finite-dimensional matrix inequalities, presented in the first theorem in this paper. Next, the second theorem is proposed to convexify the bilinear matrix inequality (BMI) in the synthesis condition of the first theorem using [31, Projection Lemma], without weakening the parameters of our KF. To further reduce conservatism, we set forth an algorithm that can compute the BMI iteratively, whose initial value can be given by a feasible solution to the synthesis condition in the second theorem. Thus, our approach eliminates the need for nonlinear SDP solvers.

The rest of the paper is organized into four sections. Section II primarily concerns the concept of KSD for dealing with the DDs in our open-loop system, and the formulation of the dissipative feedback control problem. Our main results on dissipative static controller synthesis are set out in Section III with two theorems and an iterative algorithm. Finally, the computation results of numerical examples are provided in Section IV prior to the final conclusion. To meet the page limit, some content and equations have been omitted, which can be found in the journal draft version of this article [42].

### *Notation*

Standard $p$-norm for $\mathbb{R}^n$ is defined as $\mathbb{R}^n \ni \mathbf{x} \to \|\mathbf{x}\|_p := \left(\sum_{i=1}^n |x_i|^p\right)^{1/p}$ with $p \in \mathbb{N}$. $\mathscr{M}\left(\mathcal{X};\mathbb{R}^d\right)$ stands for the set containing all measurable functions defined from Lebesgue measurable set $\mathcal{X}$ to $\mathbb{R}^d$ endowed with the Borel algebra. We use $\mathcal{C}(\mathcal{X};\mathbb{R}^n)$ to denote the Banach space of continuous functions endowed with a uniform norm $\|\boldsymbol{f}(\cdot)\|_\infty := \sup_{\tau\in\mathcal{X}} \|\boldsymbol{f}(\tau)\|_2$. We also define $\mathscr{L}^p(\mathcal{X};\mathbb{R}^n) := \{\boldsymbol{f}(\cdot) \in \mathscr{M}(\mathcal{X};\mathbb{R}^n) : \|\boldsymbol{f}(\cdot)\|_p < +\infty\}$ with $\mathcal{X} \subseteq \mathbb{R}$ and the semi-norm $\|\boldsymbol{f}(\cdot)\|_p := \left(\int_{\mathcal{X}} \|\boldsymbol{f}(x)\|_2^p \mathsf{d}x\right)^{1/p}$, and function space $\mathscr{W}^{1,2}(\mathcal{X};\mathbb{R}^n) = \{\boldsymbol{f}(\cdot) \in \mathscr{L}^2(\mathcal{X};\mathbb{R}^n) : \boldsymbol{f}'(\cdot) \in \mathscr{L}^2(\mathcal{X};\mathbb{R}^n)\}$, where $\boldsymbol{f}'(\cdot)$ is the weak derivative of $\boldsymbol{f}(\cdot)$. We utilize notation $\widetilde{\forall} x \in \mathcal{X}$, $\mathsf{P}(x)$ to indicate that property $\mathsf{P}(x)$ holds almost everywhere for $x \in \mathcal{X}$ w.r.t the Lebesgue measure. Let $\mathsf{Sy}(X) := X + X^\top$ for any square matrix. We frequently utilize $\mathbf{Col}_{i=1}^n X_i = [X_i]_{i=1}^n := \left[X_1^\top \cdots X_i^\top \cdots X_n^\top\right]^\top$ to denote a column-wise concatenation of mathematical objects, whereas $\mathbf{Row}_{i=1}^n X_i = [\![X_i]\!]_{i=1}^n = \left[X_1 \cdots X_i \cdots X_n\right]$ is the "row vector" version. Symbol $*$ is used as an abbreviation for $[*]YX = X^\top Y X$ or $X^\top Y[*] = X^\top Y X$ or $\left[\begin{smallmatrix} A & B \\ * & C \end{smallmatrix}\right] = \left[\begin{smallmatrix} A & B \\ B^\top & C \end{smallmatrix}\right]$. $\mathsf{O}_{n,m}$ stands for an $n \times m$ zero matrix that can be abbreviated as $\mathsf{O}_n$ when $n = m$, whereas $\mathbf{0}_n$ denotes an $n \times 1$ zero column vector. We use $\oplus$ to denote $X \oplus Y = \left[\begin{smallmatrix} X & \mathsf{O}_{n,q} \\ \mathsf{O}_{p,m} & Y \end{smallmatrix}\right]$ for any $X \in \mathbb{R}^{n\times m}$, $Y \in \mathbb{R}^{p\times q}$ with its $n$-ary form $\mathbf{diag}_{i=1}^\nu X_i = X_1 \oplus X_2 \oplus \cdots \oplus X_\nu$. Notation $\otimes$ stands for the Kronecker product. The order of matrix operations is defined as *matrix (scalar) multiplications* $>$ $\otimes > \oplus = \mathbf{diag} > +$. Finally, we use $[\,]$, to represent empty matrices [32, See I.7] following the same definition and rules in MATLAB©.

## II. Problem Formulation

### *A. Open-Loop LTDS*

In this paper, we deal with an LTDS of the form

$$
\begin{aligned}
\dot{\boldsymbol{x}}(t) &= \sum_{i=0}^{\nu} A_i \boldsymbol{x}(t - r_i) + \int_{-r}^{0} \widetilde{A}(\tau)\boldsymbol{x}(t+\tau)\mathsf{d}\tau \\
&\quad + \sum_{i=0}^{\nu} B_i \boldsymbol{u}(t - r_i) + \int_{-r}^{0} \widetilde{B}(\tau)\boldsymbol{u}(t+\tau)\mathsf{d}\tau + D_1\boldsymbol{w}(t), \\
\boldsymbol{z}(t) &= \sum_{i=0}^{\nu} C_i \boldsymbol{x}(t - r_i) + \int_{-r}^{0} \widetilde{C}(\tau)\boldsymbol{x}(t+\tau)\mathsf{d}\tau \\
&\quad + \sum_{i=0}^{\nu} \mathfrak{B}_i \boldsymbol{u}(t - r_i) + \int_{-r}^{0} \widetilde{\mathfrak{B}}(\tau)\boldsymbol{u}(t+\tau)\mathsf{d}\tau + D_2\boldsymbol{w}(t), \\
&\forall \theta \in [-r, 0], \quad \boldsymbol{x}(t_0 + \theta) = \boldsymbol{\psi}(\theta),
\end{aligned} \tag{1}
$$

with a quadratic supply rate function (SRF)

$$
\begin{gathered}
\mathsf{s}(\boldsymbol{z}(t), \boldsymbol{w}(t)) = \begin{bmatrix} \boldsymbol{z}(t) \\ \boldsymbol{w}(t) \end{bmatrix}^\top \begin{bmatrix} \widetilde{J}^\top J_1^{-1} \widetilde{J} & J_2 \\ * & J_3 \end{bmatrix} \begin{bmatrix} \boldsymbol{z}(t) \\ \boldsymbol{w}(t) \end{bmatrix}, \\
\widetilde{J}^\top J_1^{-1} \widetilde{J} \preceq 0, \quad J_1^{-1} \prec 0, \quad \widetilde{J} \in \mathbb{R}^{m\times m}, \\
J_2 \in \mathbb{R}^{m\times q}, \quad J_3 \in \mathbb{S}^q,
\end{gathered} \tag{2}
$$

where the functional differential equation (FDE) in (1) holds for $t \geq t_0 \in \mathbb{R}$ almost everywhere w.r.t the Lebesgue measure. The initial condition is $\boldsymbol{\psi}(\cdot) \in \mathcal{C}(\mathcal{J};\mathbb{R}^n)$ with $\mathcal{J} := [-r, 0]$, and delay values $r = r_\nu > \cdots > r_2 > r_1 > r_0 = 0$ are known with $\nu \in \mathbb{N}$. $\boldsymbol{x}(t) \in \mathbb{R}^n$ is the solution to the FDE in the sense of Carathéodory [29, page 58], $\boldsymbol{u}(t) \in \mathbb{R}^p$ is the control input, $\boldsymbol{w}(\cdot) \in \mathscr{L}^2(\mathbb{R}_{\geq t_0};\mathbb{R}^q)$ is a disturbance, and $\boldsymbol{z}(t) \in \mathbb{R}^m$ is the regulated output, where $n, m, p, q \in \mathbb{N}$. Finally, the DDs in (1) satisfy

$$
\begin{aligned}
&\widetilde{A}(\cdot) \in \mathscr{L}^2(\mathcal{J};\mathbb{R}^{n\times n}), \quad \widetilde{C}(\cdot) \in \mathscr{L}^2(\mathcal{J};\mathbb{R}^{m\times n}), \\
&\widetilde{B}(\cdot) \in \mathscr{L}^2(\mathcal{J};\mathbb{R}^{n\times p}), \quad \widetilde{\mathfrak{B}}(\cdot) \in \mathscr{L}^2(\mathcal{J};\mathbb{R}^{m\times p}).
\end{aligned} \tag{3}
$$

The integrals in (1) can always be decomposed as

$$
\begin{aligned}
\int_{-r}^{0} \widetilde{A}(\tau)\boldsymbol{x}(t+\tau)\mathsf{d}\tau &= \sum_{i=1}^{\nu} \int_{\mathcal{I}_i} \widetilde{A}_i(\tau)\boldsymbol{x}(t+\tau)\mathsf{d}\tau \\
\int_{-r}^{0} \widetilde{C}(\tau)\boldsymbol{x}(t+\tau)\mathsf{d}\tau &= \sum_{i=1}^{\nu} \int_{\mathcal{I}_i} \widetilde{C}_i(\tau)\boldsymbol{x}(t+\tau)\mathsf{d}\tau \\
\int_{-r}^{0} \widetilde{B}(\tau)\boldsymbol{u}(t+\tau)\mathsf{d}\tau &= \sum_{i=1}^{\nu} \int_{\mathcal{I}_i} \widetilde{B}_i(\tau)\boldsymbol{u}(t+\tau)\mathsf{d}\tau \\
\int_{-r}^{0} \widetilde{\mathfrak{B}}(\tau)\boldsymbol{u}(t+\tau)\mathsf{d}\tau &= \sum_{i=1}^{\nu} \int_{\mathcal{I}_i} \widetilde{\mathfrak{B}}_i(\tau)\boldsymbol{u}(t+\tau)\mathsf{d}\tau
\end{aligned} \tag{4}
$$

using $\mathcal{I}_i = [-r_i, -r_{i-1}]$ and matrix-valued functions

$$
\begin{aligned}
&\widetilde{A}_i(\cdot) \in \mathcal{L}^2(\mathcal{I}_i\,;\mathbb{R}^{n\times n}), \quad \widetilde{C}_i(\cdot) \in \mathcal{L}^2(\mathcal{I}_i\,;\mathbb{R}^{m\times n}),\\
&\widetilde{B}_i(\cdot) \in \mathcal{L}^2(\mathcal{I}_i\,;\mathbb{R}^{n\times p}), \quad \widetilde{\mathfrak{B}}_i(\cdot) \in \mathcal{L}^2(\mathcal{I}_i\,;\mathbb{R}^{m\times p}).
\end{aligned} \tag{5}
$$

for all $i \in \mathbb{N}_\nu := \{1, \ldots, \nu\}$. The structure of (1) is selected based on the general LTDSs [29] expressed as Lebesgue-Stieltjes integrals $\dot{\boldsymbol{x}}(t) = \int_{-r}^{0} \mathsf{d}\,[A(\tau)]\,\boldsymbol{x}(t+\tau) + \int_{-r}^{0} \mathsf{d}\,[B(\tau)]\,\boldsymbol{u}(t+\tau)$.

*Remark 1:* The expression of (1) is sufficiently general to describe most LTDSs from a mathematical perspective [25, Eq. (2.1)], including LTDSs with general input delays. A wide variety of cybernetic systems with general DDs can be modeled by (1) such as the characterization of event-triggered mechanisms [5], networked control systems [19], or chemical reaction networks [33, eq.(30)], etc.

We formulated $\mathsf{s}(\boldsymbol{z}(t), \boldsymbol{w}(t))$ in (2) based on the paradigm established in [34] with minor modifications. The function can describe multiple performance criteria such as

- $\mathcal{L}^2$ gain performance: $J_1 = -\gamma I_m$, $\widetilde{J} = I_m$, $J_2 = \mathsf{O}_{m,q}$, $J_3 = \gamma I_q$ with $\gamma > 0$
- Strict Passivity: $J_1 \prec 0$, $\widetilde{J} = \mathsf{O}_m$, $J_2 = I_m$, $J_3 = \mathsf{O}_m$
- Other sector constraints in [35, Table 1].

To address the challenges arising from the infinite-dimensionality of matrix-valued functions such as those in (1), we introduced the concept of KSD in [28]. The following proposition provides the first ingredient of the concept of KSD specifically for the matrix-valued functions in (1).

*Proposition 1:* (5) holds if and only if there exist $\boldsymbol{f}_i(\cdot) \in \mathcal{W}^{1,2}(\mathcal{I}_i\,;\mathbb{R}^{d_i})$, $\boldsymbol{\varphi}_i(\cdot) \in \mathcal{L}^2(\mathcal{I}_i\,;\mathbb{R}^{\delta_i})$, $\boldsymbol{\phi}_i(\cdot) \in \mathcal{L}^2(\mathcal{I}_i\,;\mathbb{R}^{\mu_i})$ and constant matrices $M_i \in \mathbb{R}^{d_i\times\varkappa_i}$, $\widehat{A}_i \in \mathbb{R}^{n\times\kappa_i n}$, $\widehat{B}_i \in \mathbb{R}^{n\times\kappa_i p}$, $\widehat{C}_i \in \mathbb{R}^{m\times\kappa_i n}$, $\widehat{\mathfrak{B}}_i \in \mathbb{R}^{m\times\kappa_i p}$ such that

$$
\widetilde{A}_i\,(\tau) = \widehat{A}_i\,(\boldsymbol{g}_i(\tau)\otimes I_n)\,, \qquad \widetilde{B}_i\,(\tau) = \widehat{B}_i\,(\boldsymbol{g}_i(\tau)\otimes I_p) \tag{6}
$$

$$
\widetilde{C}_i\,(\tau) = \widehat{C}_i\,(\boldsymbol{g}_i(\tau)\otimes I_n)\,, \qquad \widetilde{\mathfrak{B}}_i\,(\tau) = \widehat{\mathfrak{B}}_i\,(\boldsymbol{g}_i(\tau)\otimes I_p) \tag{7}
$$

$$
\frac{\mathsf{d}\boldsymbol{f}_i(\tau)}{\mathsf{d}\tau} = M_i\boldsymbol{h}_i(\tau), \qquad \boldsymbol{h}_i(\tau) = \begin{bmatrix}\boldsymbol{\varphi}_i(\tau)\\ \boldsymbol{f}_i(\tau)\end{bmatrix} \tag{8}
$$

$$
\mathsf{G}_i = \int_{\mathcal{I}_i} \boldsymbol{g}_i(\tau)\boldsymbol{g}_i^\top(\tau)\mathsf{d}\tau \succ 0, \qquad \boldsymbol{g}_i(\tau) = \begin{bmatrix}\boldsymbol{\phi}_i(\tau)\\ \boldsymbol{h}_i(\tau)\end{bmatrix} \tag{9}
$$

for all $\tau \in \mathcal{I}_i$ and $i \in \mathbb{N}_\nu$, where $\kappa_i = d_i + \delta_i + \mu_i$, $\varkappa_i = d_i + \delta_i$ with indices $d_i \in \mathbb{N}$ and $\delta_i, \mu_i \in \mathbb{N}\cup\{0\}$. The derivatives in (8) are weak derivatives [36].

*Proof:* Please see the proof of [42, Proposition 1] for a more general version of the KSD. ■

We have distinguished $\boldsymbol{\varphi}_i(\cdot)$ from $\boldsymbol{\phi}_i(\cdot)$ in $\boldsymbol{g}_i(\cdot)$, since $\boldsymbol{\phi}_i(\cdot)$ are linearly approximated by $\boldsymbol{h}_i(\cdot)$ in the following as

$$
\forall i \in \mathbb{N}_\nu,\ \forall\tau \in \mathcal{I}_i,\ \boldsymbol{\phi}_i(\tau) = \Gamma_i\,\mathfrak{H}_i^{-1}\boldsymbol{h}_i(\tau) + \boldsymbol{\varepsilon}_i(\tau) \tag{10}
$$

where $\Gamma_i := \int_{\mathcal{I}_i} \boldsymbol{\phi}_i(\tau)\boldsymbol{h}_i^\top(\tau)\mathsf{d}\tau \in \mathbb{R}^{\mu_i\times\varkappa_i}$ and $\mathfrak{H}_i := \int_{\mathcal{I}_i} \boldsymbol{h}_i(\tau)\boldsymbol{h}_i^\top(\tau)\mathsf{d}\tau$. Note that $\mathfrak{H}_i \succ 0$ is implied by (9). Similarly, we define $\boldsymbol{\varepsilon}_i(\tau) = \boldsymbol{\phi}_i(\tau) - \Gamma_i\,\mathfrak{H}_i^{-1}\boldsymbol{h}_i(\tau)$ as the approximation errors, and

$$
\begin{aligned}
\mathbb{S}^{\mu_i} \ni \mathfrak{E}_i = \int_{\mathcal{I}_i} \boldsymbol{\varepsilon}_i(\tau)\boldsymbol{\varepsilon}_i^\top(\tau)\mathsf{d}\tau = \int_{\mathcal{I}_i} \boldsymbol{\phi}_i(\tau)\boldsymbol{\phi}_i^\top(\tau)\mathsf{d}\tau\\
- \Gamma_i\,\mathfrak{H}_i^{-1}\Gamma_i^\top \succ 0
\end{aligned} \tag{11}
$$

is utilized to measure these errors, where $\mathfrak{E}_i \succ 0$ is inferred from the properties in [10, Eq. (18)]. By (10), we have

$$
\begin{aligned}
\boldsymbol{g}_i(\tau) &= \begin{bmatrix}\boldsymbol{\phi}_i(\tau)\\ \boldsymbol{h}_i(\tau)\end{bmatrix} = \begin{bmatrix}\Gamma_i\mathfrak{H}_i^{-1}\boldsymbol{h}_i(\tau)\\ \boldsymbol{h}_i(\tau)\end{bmatrix} + \begin{bmatrix}\boldsymbol{\varepsilon}_i(\tau)\\ \boldsymbol{0}_{\varkappa_i}\end{bmatrix}\\
&= \widehat{\Gamma}_i\boldsymbol{h}_i(\tau) + \widetilde{I}_i\boldsymbol{\varepsilon}_i(\tau),
\end{aligned} \tag{12}
$$

$$
\widehat{\Gamma}_i = \begin{bmatrix}\Gamma_i\mathfrak{H}_i^{-1}\\ I_{\varkappa_i}\end{bmatrix} \in \mathbb{R}^{\kappa_i\times\varkappa_i}, \quad \widetilde{I}_i = \begin{bmatrix}I_{\mu_i}\\ \mathsf{O}_{\varkappa_i,\mu_i}\end{bmatrix} \in \mathbb{R}^{\kappa_i\times\mu_i}
$$

with constant matrices $\Gamma_i, \mathfrak{H}_i$ in (10). By replacing $\boldsymbol{g}_i(\tau)$ in (6)–(7) with the terms in (12), the formulation of the KSD concept for the matrix-valued functions in (5) is completed.

*Remark 2:* We refer readers to our journal draft [42, subsection 3.1] and the tutorial in [42, Appendix E] for the motivation and explanations on the structure of KSD.

### *B. Derivation of Closed-Loop System*

Inspired by the state variable $z(t,\tau)$ in [37], we introduce $\boldsymbol{\chi}(t,\theta) = [\boldsymbol{x}(t + \acute{r}_i\theta - r_{i-1})]_{i=1}^{\nu} \in \mathbb{R}^{n\nu}$ with $\theta \in [-1, 0]$ and $\acute{r}_i = r_i - r_{i-1}$. Assuming that $\boldsymbol{x}(t)$ can be measured for feedback, we apply a static state controller $\boldsymbol{u}(t) = K\boldsymbol{x}(t)$ to (1) and utilize the decompositions in Proposition 1, where $K \in \mathbb{R}^{p\times n}$ is the gain to be computed. Then the expression of the closed-loop system (CLS) is given as

$$
\begin{aligned}
&\dot{\boldsymbol{x}}(t) = (A_0 + B_0K)\,\boldsymbol{x}(t) + \llbracket (A_i + B_iK) \rrbracket_{i=1}^{\nu}\boldsymbol{\chi}(t,-1)\\
&+ \sum_{i=1}^{\nu}\int_{\mathcal{I}_i}\Big(\widehat{A}_i + \widehat{B}_i\,(I_{\kappa_i}\otimes K)\Big)\,G_i(\tau)\boldsymbol{x}(t+\tau)\mathsf{d}\tau + D_1\boldsymbol{w}(t),\\
&\boldsymbol{z}(t) = (C_0 + \mathfrak{B}_0K)\,\boldsymbol{x}(t) + \llbracket (C_i + \mathfrak{B}_iK) \rrbracket_{i=1}^{\nu}\boldsymbol{\chi}(t,-1)\\
&+ \sum_{i=1}^{\nu}\int_{\mathcal{I}_i}\Big(\widehat{C}_i + \widehat{\mathfrak{B}}_i\,(I_{\kappa_i}\otimes K)\Big)\,G_i(\tau)\boldsymbol{x}(t+\tau)\mathsf{d}\tau + D_2\boldsymbol{w}(t),\\
&\forall\theta \in \mathcal{J},\ \ \boldsymbol{x}(t_0+\theta) = \boldsymbol{\psi}(\theta),\ \ \boldsymbol{\psi}(\cdot) \in \mathcal{C}(\mathcal{J};\mathbb{R}^n)
\end{aligned} \tag{13}
$$

where $G_i(\tau) = (\boldsymbol{g}_i(\tau)\otimes I_n)$ with $\boldsymbol{g}_i(\cdot)$ in (9), following from a use of the property of Kronecker products

$$
\begin{aligned}
\forall i \in \mathbb{N}_\nu,\ &(\boldsymbol{g}_i(\tau)\otimes I_p)\,K = (\boldsymbol{g}_i(\tau)\otimes I_p)\,(1\otimes K)\\
&= I_{\kappa_i}\boldsymbol{g}_i(\tau)\otimes KI_n = (I_{\kappa_i}\otimes K)\,(\boldsymbol{g}_i(\tau)\otimes I_n)\,.
\end{aligned} \tag{14}
$$

Given $\boldsymbol{g}_i(\tau)$ in (12) in terms of $\boldsymbol{h}_i(\tau)$ and $\boldsymbol{\varepsilon}_i(\tau)$ in light of (10), identity (14) can be further expanded as

$$
(\boldsymbol{g}_i(\tau)\otimes I_p)K = \big(\widehat{\Gamma}_i\otimes K\big)H_i(\tau) + \big(\widetilde{I}_i\otimes K\big)E_i(\tau), \tag{15}
$$

$$
\boldsymbol{g}_i(\tau)\otimes I_n = \big(\widehat{\Gamma}_i\otimes I_n\big)H_i(\tau) + \big(\widetilde{I}_i\otimes I_n\big)E_i(\tau) \tag{16}
$$

following from similar steps in [28, Eq.(18)], where $H_i(\tau) = \boldsymbol{h}_i(\tau)\otimes I_n$ and $E_i(\tau) = \boldsymbol{\varepsilon}_i(\tau)\otimes I_n$. By (15)–(16) and [28, Eq.(14a)] with matrices $\mathfrak{H}_i \succ 0$, $\mathfrak{E}_i \succ 0$ in (10) and (11), we can rewrite all the DD integral matrices in (13) as

$$
\begin{aligned}
&\Big[\widehat{A}_i + \widehat{B}_i\,(I_{\kappa_i}\otimes K)\Big]\,(\boldsymbol{g}_i(\tau)\otimes I_n)\\
&= \Big[\widehat{A}_i\,(T_i\otimes I_n) + \widehat{B}_i\,(T_i\otimes K)\Big]\,\big[\mathfrak{L}_{\boldsymbol{h}_i}^{-1}\boldsymbol{h}_i(\tau)\otimes I_n\big]\\
&\quad + \Big[\widehat{A}_i\,\Big(\widetilde{T}_i\otimes I_n\Big) + \widehat{B}_i\,\Big(\widetilde{T}_i\otimes K\Big)\Big]\,\big[\mathfrak{L}_{\boldsymbol{\varepsilon}_i}^{-1}\boldsymbol{\varepsilon}_i(\tau)\otimes I_n\big]\,,
\end{aligned} \tag{17}
$$

$$
\Big[\widehat{C}_i + \widehat{\mathfrak{B}}_i\,(I_{\kappa_i}\otimes K)\Big]\,(\boldsymbol{g}_i(\tau)\otimes I_n) \tag{18}
$$

$$
= \left[\widehat{C}_i \left(T_i \otimes I_n\right) + \widehat{\mathfrak{B}}_i \left(T_i \otimes K\right)\right] \left[\mathfrak{L}_{\boldsymbol{h}_i}^{-1} \boldsymbol{h}_i(\tau) \otimes I_n\right]
+ \left[\widehat{C}_i \left(\widetilde{T}_i \otimes I_n\right) + \widehat{\mathfrak{B}}_i \left(\widetilde{T}_i \otimes K\right)\right] \left[\mathfrak{L}_{\boldsymbol{\varepsilon}_i}^{-1} \boldsymbol{\varepsilon}_i(\tau) \otimes I_n\right]
$$

with matrices $T_i = \widehat{\Gamma}_i \mathfrak{L}_{\boldsymbol{h}_i} = \begin{bmatrix} \Gamma_i \mathfrak{L}_{\boldsymbol{h}_i}^{-\top} \\ \mathfrak{L}_{\boldsymbol{h}_i} \end{bmatrix}$ and $\widetilde{T}_i = \widetilde{I}_i \mathfrak{L}_{\boldsymbol{\varepsilon}_i}$ for all $i \in \mathbb{N}_\nu$, where $\mathfrak{L}_{\boldsymbol{h}_i} \mathfrak{L}_{\boldsymbol{h}_i}^\top = \mathfrak{H}_i$ and $\mathfrak{L}_{\boldsymbol{\varepsilon}_i} \mathfrak{L}_{\boldsymbol{\varepsilon}_i}^\top = \mathfrak{E}_i$ are the lower triangular matrices from the Cholesky decompositions. Now applying (17)–(18) and the properties in [28, Eq.(14c)] to the DD integrals in (13) further gives

$$
\begin{aligned}
&\sum_{i=1}^{\nu} \int_{\mathcal{I}_i} \left(\widehat{A}_i + \widehat{B}_i \left(I_{\kappa_i} \otimes K\right)\right) \left(\boldsymbol{g}_i(\tau) \otimes I_n\right) \boldsymbol{x}(t+\tau) \mathsf{d}\tau \\
&= \left[\!\left[\widehat{A}_i \left(T_i \otimes I_n\right) + \widehat{B}_i \left(T_i \otimes K\right)\right]\!\right]_{i=1}^{\nu} \boldsymbol{\xi}(t) \\
&\quad + \left[\!\left[\widehat{A}_i \left(\widetilde{T}_i \otimes I_n\right) + \widehat{B}_i \left(\widetilde{T}_i \otimes K\right)\right]\!\right]_{i=1}^{\nu} \boldsymbol{e}(t),
\end{aligned}
\tag{19}
$$

$$
\begin{aligned}
&\sum_{i=1}^{\nu} \int_{\mathcal{I}_i} \left(\widehat{C}_i + \widehat{\mathfrak{B}}_i \left(I_{\kappa_i} \otimes K\right)\right) \left(\boldsymbol{g}_i(\tau) \otimes I_n\right) \boldsymbol{x}(t+\tau) \mathsf{d}\tau \\
&= \left[\!\left[\widehat{C}_i \left(T_i \otimes I_n\right) + \widehat{\mathfrak{B}}_i \left(T_i \otimes K\right)\right]\!\right]_{i=1}^{\nu} \boldsymbol{\xi}(t) \\
&\quad + \left[\!\left[\widehat{C}_i \left(\widetilde{T}_i \otimes I_n\right) + \widehat{\mathfrak{B}}_i \left(\widetilde{T}_i \otimes K\right)\right]\!\right]_{i=1}^{\nu} \boldsymbol{e}(t),
\end{aligned}
\tag{20}
$$

where

$$
\begin{aligned}
\boldsymbol{\xi}(t) &= \left[\int_{\mathcal{I}_i} \left(\mathfrak{L}_{\boldsymbol{h}_i}^{-1} \boldsymbol{h}_i(\tau) \otimes I_n\right) \boldsymbol{x}(t+\tau) \mathsf{d}\tau\right]_{i=1}^{\nu}, \\
\boldsymbol{e}(t) &= \left[\int_{\mathcal{I}_i} \left(\mathfrak{L}_{\boldsymbol{\varepsilon}_i}^{-1} \boldsymbol{\varepsilon}_i(\tau) \otimes I_n\right) \boldsymbol{x}(t+\tau) \mathsf{d}\tau\right]_{i=1}^{\nu}.
\end{aligned}
\tag{21}
$$

Finally, by utilizing the identities in (19)–(21) with the properties in [28, Lemma 1] on (13), our CLS becomes

$$
\begin{aligned}
\dot{\boldsymbol{x}}(t) &= \left(\mathbf{A} + \mathbf{B}_1 \left[\left(I_\beta \otimes K\right) \oplus \mathsf{O}_q\right]\right) \boldsymbol{\vartheta}(t), \quad \forall t \geq t_0 \\
\boldsymbol{z}(t) &= \left(\mathbf{C} + \mathbf{B}_2 \left[\left(I_\beta \otimes K\right) \oplus \mathsf{O}_q\right]\right) \boldsymbol{\vartheta}(t), \\
\mathbf{x}_{t_0}(\theta) &= \boldsymbol{x}(t_0 + \theta) = \boldsymbol{\psi}(\theta), \qquad \forall \theta \in \mathcal{J}
\end{aligned}
\tag{22}
$$

with $t_0$ and $\boldsymbol{\psi}(\cdot)$ in (1), where $\beta = 1 + \nu + \kappa$ with $\kappa = \sum_{i=1}^{\nu} \kappa_i$ and $\kappa_i = d_i + \delta_i + \mu_i$ in Proposition 1, and

$$
\mathbf{A} = \left[\left[\!\left[A_i\right]\!\right]_{i=0}^{\nu} \quad \left[\!\left[\widehat{A}_i \left(T_i \otimes I_n\right)\right]\!\right]_{i=1}^{\nu} \cdots
\cdots \left[\!\left[\widehat{A}_i \left(\widetilde{T}_i \otimes I_n\right)\right]\!\right]_{i=1}^{\nu} \quad D_1\right], \tag{23}
$$

$$
\mathbf{B}_1 = \left[\left[\!\left[B_i\right]\!\right]_{i=0}^{\nu} \left[\!\left[\widehat{B}_i \left(T_i \otimes I_p\right)\right]\!\right]_{i=1}^{\nu} \cdots
\cdots \left[\!\left[\widehat{B}_i \left(\widetilde{T}_i \otimes I_p\right)\right]\!\right]_{i=1}^{\nu} \quad \mathsf{O}_{n,q}\right], \tag{24}
$$

$$
\mathbf{C} = \left[\left[\!\left[C_i\right]\!\right]_{i=0}^{\nu} \quad \left[\!\left[\widehat{C}_i \left(T_i \otimes I_n\right)\right]\!\right]_{i=1}^{\nu} \cdots
\cdots \left[\!\left[\widehat{C}_i \left(\widetilde{T}_i \otimes I_n\right)\right]\!\right]_{i=1}^{\nu} \quad D_2\right], \tag{25}
$$

$$
\mathbf{B}_2 = \left[\left[\!\left[\mathfrak{B}_i\right]\!\right]_{i=0}^{\nu} \quad \left[\!\left[\widehat{\mathfrak{B}}_i \left(T_i \otimes I_p\right)\right]\!\right]_{i=1}^{\nu} \cdots
\cdots \left[\!\left[\widehat{\mathfrak{B}}_i \left(\widetilde{T}_i \otimes I_p\right)\right]\!\right]_{i=1}^{\nu} \quad \mathsf{O}_{m,q}\right], \tag{26}
$$

$$
\boldsymbol{\omega}(t) = \begin{bmatrix} \boldsymbol{x}^\top(t) & \boldsymbol{\chi}^\top(t,-1) & \boldsymbol{\xi}^\top(t) \end{bmatrix}^\top, \tag{27}
$$

$$
\boldsymbol{\vartheta}(t) = \begin{bmatrix} \boldsymbol{\omega}^\top(t) & \boldsymbol{e}^\top(t) & \boldsymbol{w}^\top(t) \end{bmatrix}^\top. \tag{28}
$$

We define $\boldsymbol{\vartheta}(t)$ so that the terms in (13) can be expressed as the products of the matrices in (23)–(26) and $\boldsymbol{\vartheta}(t)$.

*Remark 3:* In contrast to [28], where the matrix square roots $\sqrt{\mathfrak{H}_i^{-1}}$ and $\sqrt{\mathfrak{E}_i^{-1}}$ were used to normalize $\boldsymbol{h}_i(\cdot)$ and $\boldsymbol{\varepsilon}_i(\cdot)$, we adopt Cholesky decompositions in this paper, which are considerably less expensive to compute and yield better-conditioned SDP formulations.

## III. Controller Synthesis under Dissipativity

The first theorem addressing the dissipative state feedback problem is stated below, whose proof can be found in [42, Appendix C]. The stability criteria and definition of dissipativity employed in the proof are presented in [42, Lemma 1] and [42, Definition 1].

*Theorem 1:* Let all the parameters in Proposition 1 be given. Then the CLS in (22) with SRF (2) is dissipative, and the origin of (22) with $\boldsymbol{w}(t) \equiv \mathbf{0}_q$ is globally exponentially stable if there exist $K \in \mathbb{R}^{p \times n}$, $P_1 \in \mathbb{S}^n$, $P_2 \in \mathbb{R}^{n \times \varrho}$, $P_3 \in \mathbb{S}^\varrho$ and $Q_i, R_i \in \mathbb{S}^n$, $i \in \mathbb{N}_\nu$ such that

$$
\begin{bmatrix} P_1 & P_2 \\ * & P_3 \end{bmatrix} + \left[\mathsf{O}_n \oplus \left(\operatorname*{\mathbf{diag}}_{i=1}^{\nu} I_{d_i} \otimes Q_i\right)\right] \succ 0, \tag{29}
$$

$$
\mathbf{Q} = \operatorname*{\mathbf{diag}}_{i=1}^{\nu} Q_i \succ 0, \quad \mathbf{R} = \operatorname*{\mathbf{diag}}_{i=1}^{\nu} R_i \succ 0, \tag{30}
$$

$$
\begin{bmatrix} \boldsymbol{\Psi} & \boldsymbol{\Sigma}^\top \widetilde{J}^\top \\ * & J_1 \end{bmatrix} = \mathsf{Sy}\left[\mathbf{P}^\top \boldsymbol{\Pi}\right] + \boldsymbol{\Phi} \prec 0, \tag{31}
$$

where $\varrho = dn$, $d = \sum_{i=1}^{\nu} d_i$, with $d_i = \mathsf{dim}(\boldsymbol{f}_i(\tau))$ and $\boldsymbol{\Sigma} = \mathbf{C} + \mathbf{B}_2 \left[\left(I_\beta \otimes K\right) \oplus \mathsf{O}_q\right]$ with $\mathbf{C}, \mathbf{B}_2$ in (25)–(26), and

$$
\boldsymbol{\Psi} = \mathsf{Sy}\left(S^\top \begin{bmatrix} P_1 & P_2 \\ * & P_3 \end{bmatrix} \begin{bmatrix} & \boldsymbol{\Omega} \\ \mathbf{M} \otimes I_n & \mathsf{O}_{dn,(\mu n+q)} \end{bmatrix} \cdots
- \begin{bmatrix} \mathsf{O}_{(\beta n),m} \\ J_2^\top \end{bmatrix} \boldsymbol{\Sigma}\right) + \Xi, \tag{32}
$$

$$
S = \begin{bmatrix} I_n & \mathsf{O}_{n,\nu n} & \mathsf{O}_{n,\varkappa n} & \mathsf{O}_{n,\mu n} & \mathsf{O}_{n,q} \\ \mathsf{O}_{dn,n} & \mathsf{O}_{dn,\nu n} & \widehat{I} & \mathsf{O}_{dn,\mu n} & \mathsf{O}_{dn,q} \end{bmatrix}, \tag{33}
$$

$$
\Xi = \left[\left(\mathbf{Q} + \mathbf{R}\Lambda\right) \oplus \mathsf{O}_n \oplus \mathsf{O}_{\kappa n} \oplus \mathsf{O}_q\right] - \cdots
\left[\mathsf{O}_n \oplus \mathbf{Q} \oplus \left[\operatorname*{\mathbf{diag}}_{i=1}^{\nu} I_{\varkappa_i} \otimes R_i\right] \oplus \left[\operatorname*{\mathbf{diag}}_{i=1}^{\nu} I_{\mu_i} \otimes R_i\right] \oplus J_3\right], \tag{34}
$$

$$
\widehat{I} = \left(\operatorname*{\mathbf{diag}}_{i=1}^{\nu} \mathfrak{L}_{\boldsymbol{f}_i}^{-1} \mathfrak{I}_i \, \mathfrak{L}_{\boldsymbol{h}_i}\right) \otimes I_n, \quad \mathfrak{I}_i = \begin{bmatrix} \mathsf{O}_{d_i,\delta_i} & I_{d_i} \end{bmatrix} \tag{35}
$$

$$
\Lambda = \operatorname*{\mathbf{diag}}_{i=1}^{\nu} \acute{r}_i I_n, \quad \acute{r}_i = r_i - r_{i-1}, \tag{36}
$$

$$
\mathbf{M} = \left[\operatorname*{\mathbf{diag}}_{i=1}^{\nu} \mathfrak{L}_{\boldsymbol{f}_i}^{-1} \boldsymbol{f}_i(-r_{i-1}) \quad \mathbf{0}_d \quad \mathsf{O}_{d,\varkappa}\right] - \cdots
\left[\mathbf{0}_d \quad \operatorname*{\mathbf{diag}}_{i=1}^{\nu} \mathfrak{L}_{\boldsymbol{f}_i}^{-1} \boldsymbol{f}_i(-r_i) \quad \operatorname*{\mathbf{diag}}_{i=1}^{\nu} \mathfrak{L}_{\boldsymbol{f}_i}^{-1} M_i \, \mathfrak{L}_{\boldsymbol{h}_i}\right] \tag{37}
$$

with $\varkappa = \sum_{i=1}^{\nu} \varkappa_i$, $\mu = \sum_{i=1}^{\nu} \mu_i$, $\kappa = \sum_{i=1}^{\nu} \kappa_i$ and $\kappa_i$, $\varkappa_i$, $\mu_i$, $M_i$ in Proposition 1, and $\boldsymbol{\Omega} := \mathbf{A} + \mathbf{B}_1 \left[\left(I_\beta \otimes K\right) \oplus \mathsf{O}_q\right]$

with $\mathbf{A}$, $\mathbf{B}_1$ in (23)–(24), and $\mathfrak{L}_{\boldsymbol{f}_i}$ is the Cholesky decomposition of $\mathfrak{F}_i = \int_{\mathcal{I}_i} \boldsymbol{f}_i(\tau)\boldsymbol{f}_i^\top(\tau)\mathsf{d}\tau = \mathfrak{L}_{\boldsymbol{f}_i}\mathfrak{L}_{\boldsymbol{f}_i}^\top \succ 0$. Moreover,

$$\mathbf{P} = \begin{bmatrix} P_1 & \mathsf{O}_{n,\nu n} & P_2\widehat{I} & \mathsf{O}_{n,(\mu n+q+m)} \end{bmatrix}, \ \mathbf{\Pi} = \begin{bmatrix} \mathbf{\Omega} & \mathsf{O}_{n,m} \end{bmatrix}, \tag{38}$$

$$\mathbf{\Phi} = \mathsf{Sy}\left( \begin{bmatrix} P_2 \\ \mathsf{O}_{\nu n,dn} \\ \widehat{I}^\top P_3 \\ \mathsf{O}_{(\mu n+q+m),dn} \end{bmatrix} \begin{bmatrix} \mathbf{M}\otimes I_n & \mathsf{O}_{dn,(\mu n+q+m)} \end{bmatrix} \cdots \right.$$
$$\left. + \begin{bmatrix} \mathsf{O}_{(\beta n),m} \\ -J_2^\top \\ \widetilde{J} \end{bmatrix} \begin{bmatrix} \mathbf{\Sigma} & \mathsf{O}_m \end{bmatrix} \right) + \Xi(\mathbf{Q},\mathbf{R}) \oplus J_1.$$

The inequality in (31) is bilinear (nonconvex) due to the products of $K$ and $P_1, P_2$. The subsequent theorem, whose proof is found in [42, Appendix E], addresses this issue by decoupling the BMI in (31) using the Projection Lemma [31], [38]. The core strategy of these methods [39], [40] is to construct convex SDP via the introduction of slack variables, while preserving the structural integrity of $P_2 \in \mathbb{R}^{n\times dn}$.

*Theorem 2:* Given $\{\alpha_i\}_{i=1}^{\beta} \subset \mathbb{R}$ and the functions and parameters in Proposition 1, then CLS (22) with SRF (2) is dissipative and the origin of (22) with $\boldsymbol{w}(t) \equiv \mathbf{0}_q$ is globally exponentially stable if there exist $X, \acute{P}_1 \in \mathbb{S}^n$, $\acute{P}_2 \in \mathbb{R}^{n\times\varrho}$, $\acute{P}_3 \in \mathbb{S}^{\varrho}$ and $\acute{Q}_i; \acute{R}_i \in \mathbb{S}^n$, $\varrho = nd$ and $V \in \mathbb{R}^{p\times n}$ such that

$$\begin{bmatrix} \acute{P}_1 & \acute{P}_2 \\ * & \acute{P}_3 \end{bmatrix} + \left[ \mathsf{O}_n \oplus \left( \operatorname*{\mathbf{diag}}_{i=1}^{\nu} I_{d_i} \otimes \acute{Q}_i \right) \right] \succ 0, \tag{39}$$

$$\acute{\mathbf{Q}} = \operatorname*{\mathbf{diag}}_{i=1}^{\nu} \acute{Q}_i \succ 0, \quad \acute{\mathbf{R}} = \operatorname*{\mathbf{diag}}_{i=1}^{\nu} \acute{R}_i \succ 0 \tag{40}$$

$$\mathsf{Sy}\left( \begin{bmatrix} I_n \\ [\alpha_i I_n]_{i=1}^{\beta} \\ \mathsf{O}_{(q+m),n} \end{bmatrix} \begin{bmatrix} -X & \acute{\mathbf{\Pi}} \end{bmatrix} \right) + \begin{bmatrix} \mathsf{O}_n & \acute{\mathbf{P}} \\ * & \acute{\mathbf{\Phi}} \end{bmatrix} \prec 0 \tag{41}$$

where $\acute{\mathbf{P}} = \begin{bmatrix} \acute{P}_1 & \mathsf{O}_{n,\nu n} & \acute{P}_2\widehat{I} & \mathsf{O}_{n,(n\mu+q+m)} \end{bmatrix}$,

$$\acute{\mathbf{\Pi}} = \begin{bmatrix} \mathbf{A}\left[(I_\beta\otimes X)\oplus I_q\right] + \mathbf{B}_1\left[(I_\beta\otimes V)\oplus \mathsf{O}_q\right] & \mathsf{O}_{n,m} \end{bmatrix}$$

with $\widehat{I}$ in (35), and matrices $\acute{\mathbf{\Phi}} = \mathbf{\Phi}(\acute{P}_2, \acute{P}_3, \acute{\mathbf{Q}}, \acute{\mathbf{R}}, \acute{\mathbf{\Sigma}})$ and $\mathbf{M}$ in (37) and $\acute{\mathbf{\Sigma}} = \mathbf{C}\left[(I_\beta\otimes X)\oplus I_q\right] + \mathbf{B}_2\left[(I_\beta\otimes V)\oplus \mathsf{O}_q\right]$ with the parameters $\mathbf{A}, \mathbf{B}_1, \mathbf{B}_2, \mathbf{C}$ in (23)–(26). Controller gain $K$ is calculated via $K = VX^{-1}$.

*Remark 4:* A qualitative analysis of the conservatism of the proposed theorems, demonstrating that the Krasovskiĭ functional in our proof parameterizes the complete Krasovskiĭ functional via basis functions $\boldsymbol{f}_i(\cdot)$, is provided in [42, Appendix D].

### A. Inner Convex Approximation of BMI

Although the constraints in Theorem 2 are convex, the introduction of $[\alpha_i I_n]_{i=1}^{\beta}$ may introduce conservatism relative to the original synthesis condition in Theorem 1. Thus, methods that can directly solve (31) are preferred. Here, we propose an offline sequential convex SDP algorithm based on the inner convex approximation strategy developed by [41].

By making use of procedures in the journal draft [42, Eq.(43)-(46)], we can conclude that (31) is implied by

$$\begin{bmatrix} \widehat{\mathbf{\Phi}} + \mathsf{Sy}\left(\widetilde{\mathbf{P}}^\top\mathbf{N} + \mathbf{P}^\top\widetilde{\mathbf{N}} - \widetilde{\mathbf{P}}^\top\widetilde{\mathbf{N}}\right) & \mathbf{P}^\top - \widetilde{\mathbf{P}}^\top & \mathbf{N}^\top - \widetilde{\mathbf{N}}^\top \\ * & -Z & \mathsf{O}_n \\ * & * & Z - I_n \end{bmatrix} \prec 0 \tag{42}$$

where $\widetilde{\mathbf{P}} = \begin{bmatrix} \widetilde{P}_1 & \mathsf{O}_{n,\nu n} & \widetilde{P}_2\widehat{I} & \mathsf{O}_{n,(n\mu+q+m)} \end{bmatrix}$,

$$\widetilde{P}_1 \in \mathbb{S}^n,\ \widetilde{P}_2 \in \mathbb{R}^{n\times\varrho},\ \mathbf{Y} = \begin{bmatrix} P_1 & P_2 \end{bmatrix},\ \widetilde{\mathbf{Y}} = \begin{bmatrix} \widetilde{P}_1 & \widetilde{P}_2 \end{bmatrix},$$
$$\mathbf{N} = \begin{bmatrix} \mathbf{B}_1 & \mathsf{O}_{n,m} \end{bmatrix}\left[(I_\beta\otimes K)\oplus \mathsf{O}_{q+m}\right],\quad Z \in \mathbb{S}^n \tag{43}$$
$$\widetilde{\mathbf{N}} = \begin{bmatrix} \mathbf{B}_1 & \mathsf{O}_{n,m} \end{bmatrix}\left[\left(I_\beta\otimes \widetilde{K}\right)\oplus \mathsf{O}_{q+m}\right],$$
$$\widehat{\mathbf{\Phi}} := \mathsf{Sy}\left(\mathbf{P}^\top\begin{bmatrix} \mathbf{A} & \mathsf{O}_{n,m} \end{bmatrix}\right) + \mathbf{\Phi}, \tag{44}$$

where (42) is a **convex** constraint if $\widetilde{\mathbf{P}}$ and $\widetilde{K}$ are known.

By following the procedures in [42, Eq.(43)-(47)] according to the expositions in [41], Algorithm 1 is established, where $f(\mathbf{x})$ is a convex objective function with $\mathbf{x}$ comprising all the decision variables in (42), and $\|\mathbf{x}\|_\infty := \max_i |x_i|$ and $\|\cdot\|_F$ is the Frobenius norm. Scalars $\rho_1, \rho_2$ and $\varepsilon$ are given constants for regularization and regulating error tolerance, respectively. Initialized from a feasible solution to the convex constraints in (39)–(41), the algorithm generates a monotone sequence and converges to a KKT point (or stationary point) based on [41, Theorem 3.2] for general nonconvex matrix inequalities, where each step or iteration solves a **convex** SDP program, thereby integrating the strengths of both Theorem 1 and Theorem 2.

**Algorithm 1:** An iterative approach for Theorem 1

**begin**
- **choose** $\{\alpha_i\}_{i=1}^{\beta} \subset \mathbb{R}$, **solve** (39)–(41) **return** $K$;
- **solve** (29)–(31) with $K$ **return** $\mathbf{Y} = \begin{bmatrix} P_1 & P_2 \end{bmatrix}$;
- **solve** (29)–(31) again with $\mathbf{Y}$ **return** $K$;
- **update** $\widetilde{\mathbf{Y}} \longleftarrow \mathbf{Y} = \begin{bmatrix} P_1 & P_2 \end{bmatrix}$, $\widetilde{K} \longleftarrow K$;
- **solve** $\min\limits_{\mathbf{x}} f(\mathbf{x}) + \rho_1\left\|\mathbf{Y} - \widetilde{\mathbf{Y}}\right\|_F^2 + \rho_2\left\|K - \widetilde{K}\right\|_F^2$ subject to (29)–(30), (42) with (43) and the parameters in Theorem 1, **return** $\mathbf{Y}$ and $K$;
- **while** $\dfrac{\left\|\begin{bmatrix} \mathbf{vec}(\mathbf{Y}) \\ \mathbf{vec}(K) \end{bmatrix} - \begin{bmatrix} \mathbf{vec}(\widetilde{\mathbf{Y}}) \\ \mathbf{vec}(\widetilde{K}) \end{bmatrix}\right\|_\infty}{\left\|\begin{bmatrix} \mathbf{vec}(\widetilde{\mathbf{Y}}) \\ \mathbf{vec}(\widetilde{K}) \end{bmatrix}\right\|_\infty + 1} \geq \varepsilon$ **do**
  - **update** $\widetilde{\mathbf{Y}} \longleftarrow \mathbf{Y}$, $\widetilde{K} \longleftarrow K$;
  - **solve** again the SDP optimization problem in the previous step, **return** $\mathbf{Y}$ and $K$;
- **end**

**end**

## IV. Numerical Experiments and Simulations

We conducted numerical experiments to demonstrate the effectiveness of the proposed framework and the advantages of the Carathéodory formulation for modeling FDEs whose parameters are subject to noise and glitches. All computations were performed using the MATLAB© platform with the package Yalmip [43] as the optimization parser, and with Mosek [44], SDPT3 [45] employed as the numerical solvers for SDP problems.

Consider a system of the form (1) with $r_1 = 1$, $r_2 = 1.7$ and the state space matrices

$$A_0 = \begin{bmatrix} -2 & 0 \\ 2 & 0.01 \end{bmatrix},\ A_1 = \begin{bmatrix} -1 & 0.1 \\ 0.2 & 0 \end{bmatrix},\ A_2 = \begin{bmatrix} -0.1 & 0 \\ 0 & -0.2 \end{bmatrix},$$

$$B_0 = \begin{bmatrix} 0 \\ 1 \end{bmatrix},\ B_1 = \begin{bmatrix} 0.01 \\ 0.1 \end{bmatrix},\ B_2 = -\begin{bmatrix} 0.1 \\ 0.1 \end{bmatrix}$$

$$\widetilde{A}_1(\tau) = \begin{bmatrix} 0.1{+}3\sin(20\tau) & 0.8\exp(\sin 20\tau){-}0.3\exp(\cos 20\tau) \\ 0.3{+}\frac{1}{\sin^2(1.2\tau)+1.0} & 3\sin(20\tau) \end{bmatrix},$$

$$\widetilde{A}_2(\tau) = \begin{bmatrix} -10\cos(18\tau) & 0.3\exp(\cos 18\tau){-}\frac{1}{\cos^2 0.7\tau+1} \\ 0.1\exp(\sin 18\tau) & 0.2{-}10\cos(18\tau) \end{bmatrix},$$

$$\widetilde{B}_1(\tau) = \begin{bmatrix} 0.01\tau{-}\frac{0.01}{\sin^2(1.2\tau)+1}{+}0.1 \\ 0.1\tau{+}\frac{0.02}{\sin^2(1.2\tau)+1} \end{bmatrix} + \varsigma(t), \tag{45}$$

$$\widetilde{B}_2(\tau) = \begin{bmatrix} 0.2\exp(\cos 18\tau){+}0.01\exp(\sin 18\tau){+}\frac{0.01}{\cos^2(0.7\tau)+1} \\ 0.1\exp(\cos 18\tau){+}0.02\exp(\sin(18\tau)) \end{bmatrix}$$

$$C_0 = \begin{bmatrix} -0.1 & 0.2 \\ 0 & 0.1 \end{bmatrix},\ C_1 = \begin{bmatrix} -0.1 & 0 \\ 0 & 0.2 \end{bmatrix},\ C_2 = \begin{bmatrix} 0 & 0.1 \\ -0.1 & 0 \end{bmatrix},$$

$$\widetilde{C}_1(\tau) = \begin{bmatrix} 0.7{+}\cos(20\tau) & \frac{1}{\sin^2 1.2\tau+1}{-}0.2 \\ 0.4{-}0.5\exp(\sin 20\tau) & 0.8{-}\sin(20\tau) \end{bmatrix},$$

$$\widetilde{C}_2(\tau) = \begin{bmatrix} 0.2{+}\sin(18\tau) & 0.3{+}\exp(\cos 18\tau) \\ 0 & 0.1{-}\frac{1}{\cos^2 0.7\tau+1} \end{bmatrix}, \mathfrak{B}_0 = \begin{bmatrix} 0 \\ 1 \end{bmatrix}$$

$$\widetilde{\mathfrak{B}}_1(\tau) = \begin{bmatrix} 0.01\tau{+}0.1\exp(\sin 20\tau){-}\frac{0.1}{\sin^2(1.2\tau)+1} \\ 0.2\exp(\sin 20\tau) \end{bmatrix},$$

$$\widetilde{\mathfrak{B}}_2(\tau) = \begin{bmatrix} 0.2\exp(\cos 18\tau){+}0.01\exp(\sin 18\tau){+}\frac{0.1}{\cos^2(0.7\tau)+1} \\ 0.02\exp(\sin 18\tau){+}\frac{0.2}{\cos^2(0.7\tau)+1} \end{bmatrix}$$

$$\mathfrak{B}_1 = \begin{bmatrix} 0.01 \\ 0.01 \end{bmatrix},\ \ \mathfrak{B}_2 = -\begin{bmatrix} 0.01 \\ 0.1 \end{bmatrix},\ D_1 = \begin{bmatrix} 0.2 \\ 0.3 \end{bmatrix},\ D_2 = \begin{bmatrix} 0.12 \\ 0.1 \end{bmatrix}$$

with $n = m = 2$, $p = q = 1$, where $\varsigma(t) = \mathbf{0}_2$ holds almost everywhere w.r.t the Lebesgue measure. Time-varying signal $\varsigma(t)$ could represent glitches or other anomalies in the input gain matrix $\widetilde{B}_1(\tau)$ and is effectively treated as zero, which could not be characterized if we were to use the traditional derivative for $\dot{\boldsymbol{x}}(t)$ in (1). This serves as an example to demonstrate the benefits of utilizing the Carathéodory framework in the modeling of LTDSs. Employing the spectral method from [46], we find that the nominal system with $\boldsymbol{w}(t) \equiv \mathbf{0}_q$ is unstable. Moreover, we employ the $\mathscr{L}^2$ gain

$$\gamma > 0,\ J_1 = -\gamma I_2,\ \widetilde{J} = I_2,\ J_2 = \mathbf{0}_2,\ J_3 = \gamma \tag{46}$$

as the performance objective $f(\mathbf{x}) = \gamma$ for the supply rate function in (2) with $\gamma$ to be minimized.

Assuming that all the system's states are measurable, our goal is to determine the controller gain of $K$ to stabilize the open-loop system in (1), while minimizing the $\mathscr{L}^2$ gain. By examining the DD kernels in (45), let

$$\boldsymbol{\phi}_1(\tau) = \begin{bmatrix} \exp(\sin 20\tau) \\ \exp(\cos 20\tau) \end{bmatrix},\quad \boldsymbol{\phi}_2(\tau) = \begin{bmatrix} \exp(\sin 18\tau) \\ \exp(\cos 18\tau) \end{bmatrix} \tag{47a}$$

$$\boldsymbol{\varphi}_1(\tau) = \frac{1}{\sin^2 1.2\tau + 1},\quad \boldsymbol{\varphi}_2(\tau) = \frac{1}{\cos^2 0.7\tau + 1} \tag{47b}$$

$$\boldsymbol{f}_1(\tau) = \begin{bmatrix} [\tau^i]_{i=0}^{\sigma_1} \\ [\sin 20i\tau]_{i=1}^{\lambda_1} \\ [\cos 20i\tau]_{i=1}^{\lambda_1} \end{bmatrix},\quad \boldsymbol{f}_2(\tau) = \begin{bmatrix} [\tau^i]_{i=0}^{\sigma_2} \\ [\sin 18i\tau]_{i=1}^{\lambda_2} \\ [\cos 18i\tau]_{i=1}^{\lambda_2} \end{bmatrix} \tag{47c}$$

for (6)–(9) with $d_i = 2\lambda_i + \sigma_i + 1$, $\mu_i = 2$, $\delta_i = 1$ and

$$M_1 = \begin{bmatrix} \mathbf{0}_{d_1} & \begin{bmatrix} \mathbf{0}_{\sigma_1}^\top & 0 \\ \mathbf{diag}_{i=1}^{\sigma_1}\, i & \mathbf{0}_{\sigma_1} \end{bmatrix} \oplus \begin{bmatrix} \mathsf{O}_{\lambda_1} & \mathbf{diag}_{i=1}^{\lambda_1}\, 20i \\ -\mathbf{diag}_{i=1}^{\lambda_1}\, 20i & \mathsf{O}_{\lambda_1} \end{bmatrix} \end{bmatrix}$$

$$M_2 = \begin{bmatrix} \mathbf{0}_{d_2} & \begin{bmatrix} \mathbf{0}_{\sigma_2}^\top & 0 \\ \mathbf{diag}_{i=1}^{\sigma_2}\, i & \mathbf{0}_{\sigma_2} \end{bmatrix} \oplus \begin{bmatrix} \mathsf{O}_{\lambda_2} & \mathbf{diag}_{i=1}^{\lambda_2}\, 18i \\ -\mathbf{diag}_{i=1}^{\lambda_2}\, 18i & \mathsf{O}_{\lambda_2} \end{bmatrix} \end{bmatrix}$$

for the relations in (8). As a result, we can construct

$$\begin{aligned}
\widehat{A}_1 &= \begin{bmatrix} 0 & 0.8 & 0 & -0.3 & 0 & 0 & 0.1 & 0 & \mathbf{0}_{2\sigma_1}^\top & 3 & 0 & \mathbf{0}_{4\lambda_1-2}^\top \\ 0 & 0 & 0 & 0 & 1 & 0 & 0.3 & 0 & \mathbf{0}_{2\sigma_1}^\top & 0 & 3 & \mathbf{0}_{4\lambda_1-2}^\top \end{bmatrix} \\
\widehat{A}_2 &= \begin{bmatrix} 0 & 0 & 0 & 0.3 & 0 & -1 & 0 & 0 & \mathbf{0}_{2\sigma_2+2\lambda_2}^\top & -10 & 0 & \mathbf{0}_{2\lambda_2-2}^\top \\ 0.1 & 0 & 0 & 0 & 0 & 0 & 0 & 0.2 & \mathbf{0}_{2\sigma_2+2\lambda_2}^\top & 0 & -10 & \mathbf{0}_{2\lambda_2-2}^\top \end{bmatrix} \\
\widehat{B}_1 &= \begin{bmatrix} 0 & 0 & -0.01 & 0.1 & 0.01 & \mathbf{0}_{\sigma_1-1+2\lambda_1}^\top \\ 0 & 0 & 0.02 & 0 & 0.1 & \mathbf{0}_{\sigma_1-1+2\lambda_1}^\top \end{bmatrix} \\
\widehat{B}_2 &= \begin{bmatrix} 0.01 & 0.2 & 0.01 & \mathbf{0}_{\sigma_2+1+2\lambda_2}^\top \\ 0.02 & 0.1 & 0 & \mathbf{0}_{\sigma_2+1+2\lambda_2}^\top \end{bmatrix} \\
\widehat{C}_1 &= \begin{bmatrix} 0 & 0 & 0 & 0 & 0 & 1 & 0.7 & -0.2 & \mathbf{0}_{2\sigma_1}^\top & 0 & 0 & \mathbf{0}_{2\lambda_1-2}^\top & 1 & 0 & \mathbf{0}_{2\lambda_1-2}^\top \\ -0.5 & 0 & 0 & 0 & 0 & 0 & 0.4 & 0.8 & \mathbf{0}_{2\sigma_1}^\top & 0 & -1 & \mathbf{0}_{2\lambda_1-2}^\top & 0 & 0 & \mathbf{0}_{2\lambda_1-2}^\top \end{bmatrix} \\
\widehat{C}_2 &= \begin{bmatrix} 0 & 0 & 0 & 1 & 0 & 0 & 0.2 & 0.3 & \mathbf{0}_{2\sigma_2}^\top & 1 & 0 & \mathbf{0}_{4\lambda_2-2}^\top \\ 0 & 0 & 0 & 0 & 0 & -1 & 0 & 0.1 & \mathbf{0}_{2\sigma_2}^\top & 0 & 0 & \mathbf{0}_{4\lambda_2-2}^\top \end{bmatrix} \\
\widehat{\mathfrak{B}}_1 &= \begin{bmatrix} 0.1 & 0 & -0.1 & 0 & 0.01 & \mathbf{0}_{\sigma_1+2\lambda_1-1}^\top \\ 0.2 & 0 & 0 & 0 & 0 & \mathbf{0}_{\sigma_1+2\lambda_1-1}^\top \end{bmatrix} \\
\widehat{\mathfrak{B}}_2 &= \begin{bmatrix} 0.01 & 0.2 & 0.1 & 0 & 0 & \mathbf{0}_{\sigma_2+2\lambda_2-1}^\top \\ 0.02 & 0 & 0.2 & 0 & 0 & \mathbf{0}_{\sigma_2+2\lambda_2-1}^\top \end{bmatrix}
\end{aligned} \tag{48}$$

to satisfy the conditions in (6)–(9).

To compute $K$, apply Theorem 2 to (22) with $\sigma_1 = \sigma_2 = \lambda_1 = \lambda_2 = 1$ and $\alpha_i = 0$, $i = 2, \ldots, \beta$, $\alpha_1 = 5$ and the parameters in (45)–(48), where $\Gamma_i$, $\mathfrak{F}_i$, $\mathfrak{H}_i$, $\mathfrak{E}_i$ are computed via the MATLAB© function `vpaintegral` that performs numerical integrations with high-level variable precision. Our SDP program yields numerical results $K = -\begin{bmatrix} 1.3794 & 1.8668 \end{bmatrix}$ with $\min \gamma = 0.8986$. This $K$ is then used for initializing Algorithm 1. After running Algorithm 1 with different numbers of iterations (NoIs) for the same system with different sets of $\sigma_i, \lambda_i$, the numerical results are listed in Tables I and II, where the spectral abscissae (SA) of the resulting CLSs with $\boldsymbol{w}(t) \equiv \mathbf{0}_q$ were calculated by the method in [46]. Our results shows that increasing $\mathsf{dim}(\boldsymbol{f}_i(\tau)) = d_i$ by stacking more functions (larger $\lambda_1, \lambda_2$) in $\boldsymbol{f}_i(\cdot) \in \mathscr{W}^{1,2}(\mathcal{I}_i\,; \mathbb{R}^{d_i})$ satisfying (9) may increase the feasibility of the synthesis conditions, leading to smaller $\min \gamma$. Moreover, it confirms that using Algorithm 1 can produce controller gains of significantly better performance than those produced by employing Theorem 2 alone ($\min \gamma = 0.8986$). This illustrates the contribution of Algorithm 1.

For numerical simulations, we consider the CLS in (22) with the parameters in (45) and controller gain $K = \begin{bmatrix} -1.5810 & -1.9805 \end{bmatrix}$ in Table II that guarantees $\min \gamma = 0.6361$. Let $t_0 = 0$ and $\boldsymbol{\psi}(\tau) = \begin{bmatrix} 5 & 3 \end{bmatrix}^\top, \tau \in [-r_2, 0]$ as the initial condition, and $w(t) = 5\sin 3\pi t(\mathit{1}(t) - \mathit{1}(t-10))$ as the disturbance where $\mathit{1}(t)$ is the Heaviside step function.

TABLE I: Controller gains $K$ with $\min\gamma$ produced with $\sigma_1 = \sigma_2 = \lambda_1 = \lambda_2 = 1,\ d_1 = d_2 = 4$.

| $K$ | $\begin{bmatrix} -1.5456 \\ -1.9359 \end{bmatrix}^\top$ | $\begin{bmatrix} -1.5365 \\ -1.9539 \end{bmatrix}^\top$ | $\begin{bmatrix} -1.5180 \\ -1.9696 \end{bmatrix}^\top$ | $\begin{bmatrix} -1.5033 \\ -1.9815 \end{bmatrix}^\top$ |
|---|---|---|---|---|
| $\min\gamma$ | 0.6573 | 0.6542 | 0.6523 | 0.6509 |
| SA | $-0.7223$ | $-0.7214$ | $-0.7224$ | $-0.7233$ |
| NoIs | 5 | 10 | 15 | 20 |

TABLE II: Controller gains $K$ with $\min\gamma$ produced with $\sigma_1 = \sigma_2 = 1,\ \lambda_1 = \lambda_2 = 2,\ d_1 = d_2 = 6$.

| $K$ | $\begin{bmatrix} -1.5538 \\ -1.9566 \end{bmatrix}^\top$ | $\begin{bmatrix} -1.5848 \\ -1.9638 \end{bmatrix}^\top$ | $\begin{bmatrix} -1.5870 \\ -1.9721 \end{bmatrix}^\top$ | $\begin{bmatrix} -1.5810 \\ -1.9805 \end{bmatrix}^\top$ |
|---|---|---|---|---|
| $\min\gamma$ | 0.6443 | 0.6398 | 0.6376 | 0.6361 |
| SA | $-0.7179$ | $-0.7113$ | $-0.7099$ | $-0.7099$ |
| NoIs | 5 | 10 | 15 | 20 |

For noise and glitch signal $\varsigma(t)$, we employ the Band-Limited White Noise block in Simulink to generate a white noise signal $n(t)$ for $\varsigma(t) = n(t)\begin{bmatrix} 1 & 1 \end{bmatrix}^\top$ in (45), where `Sample time` $= 0.002$s and the default values for `Seed` and `Noise power` were adopted. Since $n(t)$ can only be realized as a discrete sequence $n(t) = n(kT)$ within a numerical simulation environment, function $\varsigma(t) = \varsigma(kT)$ has only a **finite** number of nonzero values, which satisfies $\widetilde{\forall} t \geq 0,\ \varsigma(t) = \varsigma(kT) = 0$ and is in line with the definitions in (45). Finally, the computations were performed via the ODE solver `ode8` a fixed step size of $0.002$s.

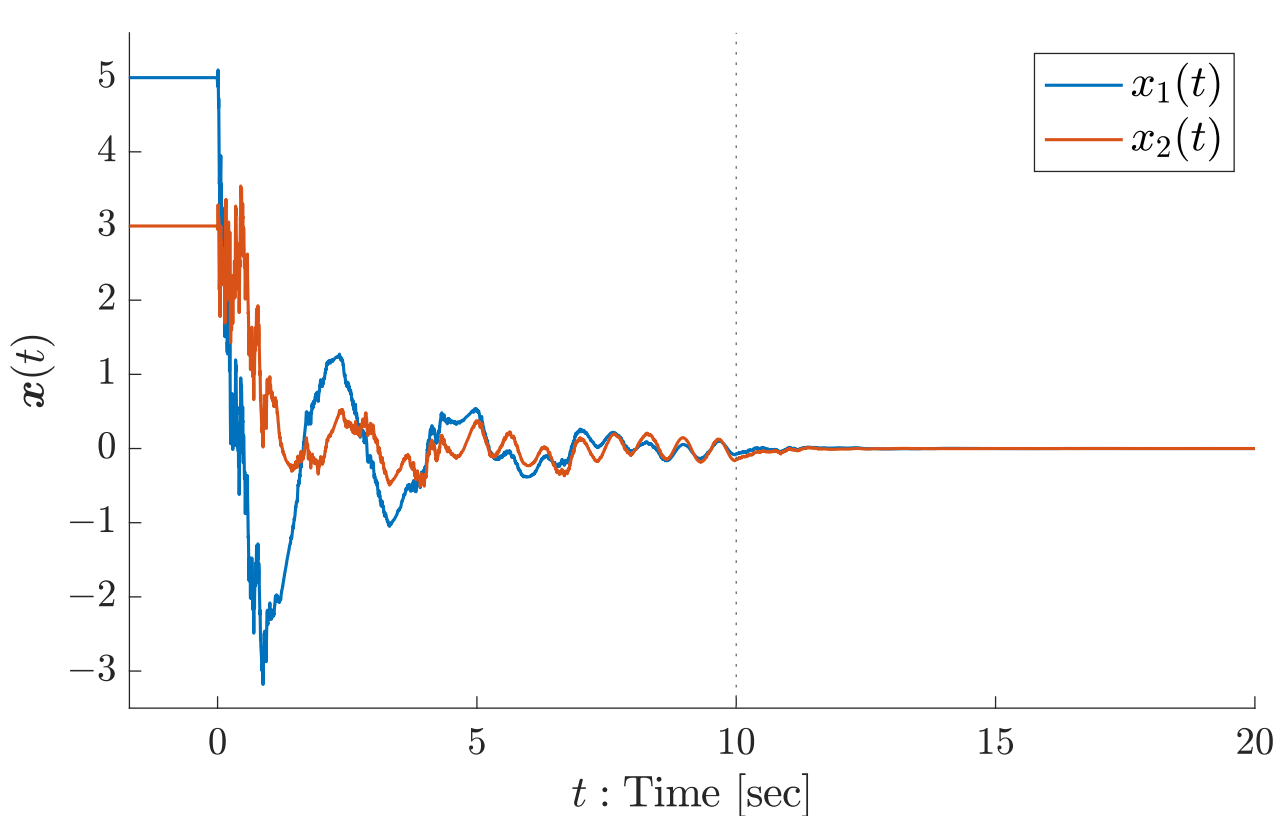

Fig. 1: Plot of $\boldsymbol{x}(t)$ with $K$ in Table II ensuring $\min\gamma = 0.6361$.

## V. Conclusion

We have proposed an SDP framework based on the KSD concept [28] for dissipative controller design of the general LTDS in (1), whose structure closely mirrors LTDSs expressed via Lebesgue-Stieltjes integrals with FDEs understood in the extended sense. The KSD overcomes the infinite dimensionality of DD-matrix kernels, yielding finite-dimensional synthesis conditions in Theorem 1 and the convex constraints in Theorem 2, while Algorithm 1 tackles the non-convexity arising from the BMI in (31). Numerical experiments confirm that the framework can effectively compute dissipative controllers for small- to medium-size systems with intricate delay structures, even when the kernel functions exhibit vastly different properties.

The proposed methodology offers a foundation for addressing delay-related practical systems such as neural networks [47], event-triggered systems [5], and other applications [48], as evidenced by [2], [5], [47], [49], [50], which utilized the decomposition in [9]—a special instance of the KSD. Moreover, the framework can serve as a blueprint for other semi-open problems in LTDSs, such as dissipative dynamical output feedback control and filtering.

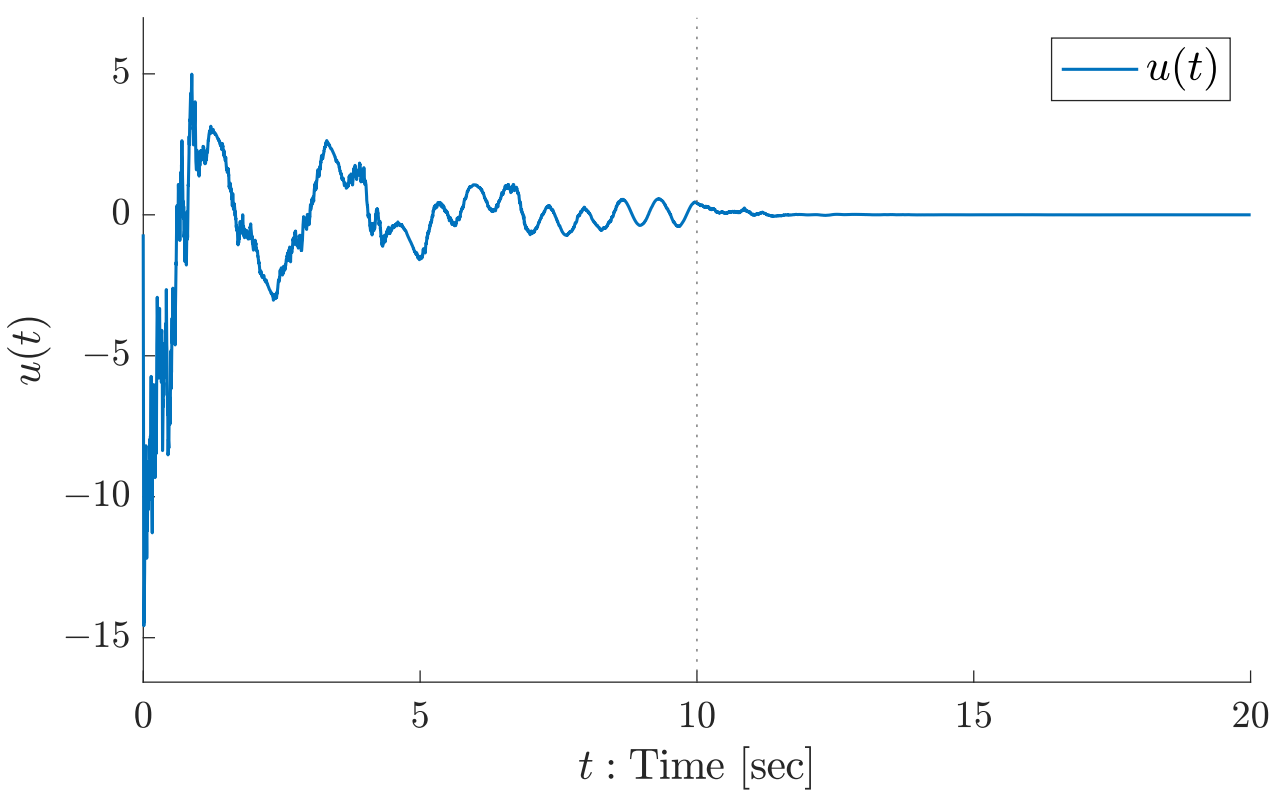

Fig. 2: Plot of $u(t) = K\boldsymbol{x}(t)$ ensuring $\min\gamma = 0.6361$.

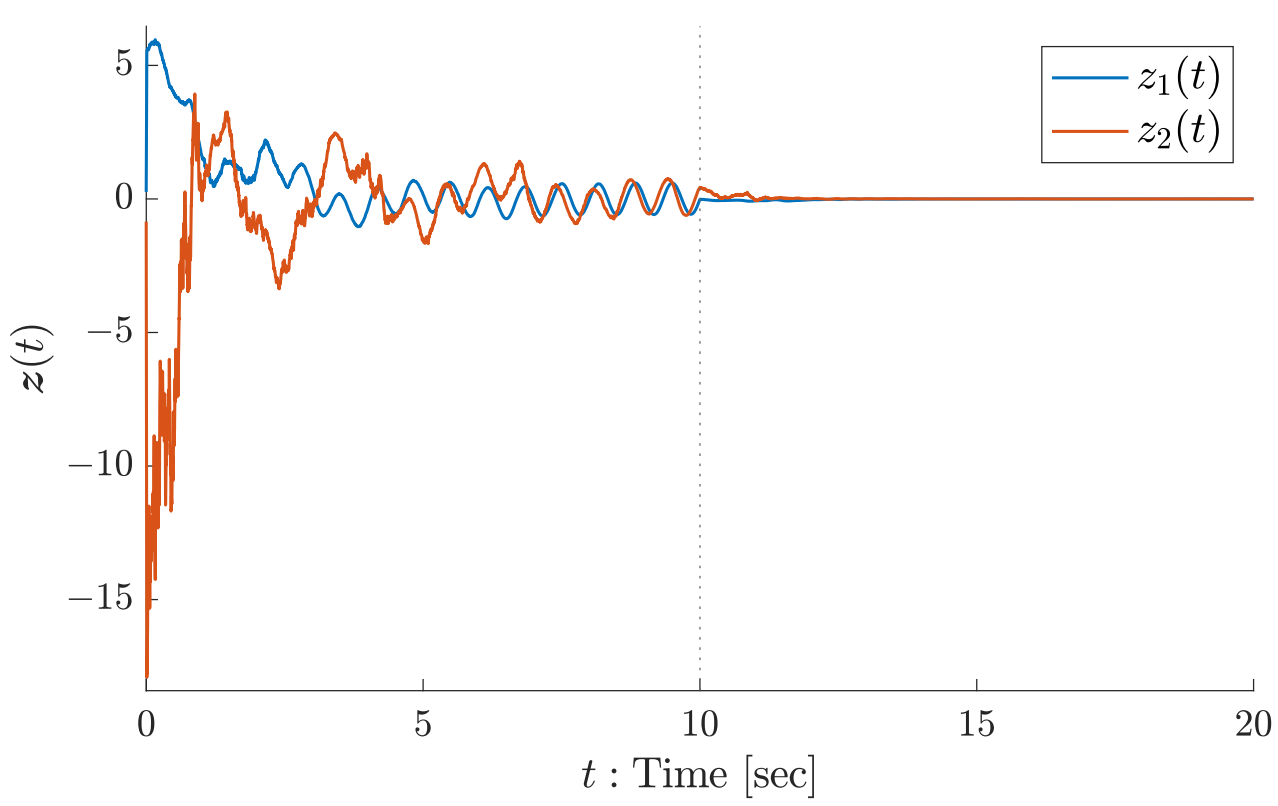

Fig. 3: Plot of $\boldsymbol{z}(t)$ in (22) ensuring $\min\gamma = 0.6361$.